\documentclass[letterpaper, 10 pt, conference]{ieeeconf}  

\usepackage{amsmath}
\usepackage{amsfonts}
\usepackage{xcolor}
\newtheorem{theorem}{Theorem}[section]
\newtheorem{lemma}[theorem]{Lemma}
\newtheorem{proposition}[theorem]{Proposition}
\newtheorem{conjecture}[theorem]{Conjecture}
\newtheorem{example}[theorem]{Example}
\newtheorem{definition}[theorem]{Definition}
\newtheorem{remark}[theorem]{Remark}

\providecommand{\noopsort}[1]{} 

\newcommand{\ext}{\text{ext}}

\DeclareMathOperator{\ran}{ran}

\usepackage{xpatch}

\xpatchcmd{\appendix}
  {Appendix \\*}
  {Appendix:\space}
  {}
  {\PackageWarning{appendix-format}{Patch failed}}

\IEEEoverridecommandlockouts                              

\title{\LARGE \bf
Controllability with and without Dissipation -- a Port-Hamiltonian Approach*
}

\author{Jochen Glück$^{1}$, Anthony Hastir$^{2}$, Birgit Jacob$^{1}$ and Hannes Wagener$^{1}$ 
\thanks{*This research was conducted with the financial support of the Deutsche Forschungsgemeinschaft (DFG, German Research Foundation). A.~Hastir has been supported by the DFG - Project-ID 532208976. All the authors acknowledge funding by the DFG, Germany
- Project-ID 531152215 - CRC 1701.}
\thanks{$^{1}$J.~Glück, B.~Jacob and H.~Wagener are with the School of Mathematics and Natural Sciences, Bergische Universit\"at Wuppertal, Gau\ss\-stra\ss e 20,
42119 Wuppertal, Germany,
        {\tt glueck@uni-wuppertal.de}, {\tt bjacob@uni-wuppertal.de}, {\tt wagener@uni-wuppertal.de}}%
\thanks{$^{2}$A.~Hastir is with the Department of Mathematics and naXys,
        University of Namur, 5000 Namur, Belgium,
        {\tt anthony.hastir@unamur.be}}%
}

\DeclareMathOperator{\spn}{span}
\DeclareMathOperator{\rank}{rk}

\begin{document}

\maketitle
\thispagestyle{plain}
\pagestyle{plain}

\begin{abstract}
    A port-Hamiltonian system with energy dissipation can be modelled by connecting a larger, non-dissipative system with a so-called closure relation. We study whether approximative controllability of the dissipative system is characterised by the same property of the larger system. The answer is negative in general, but positive for sufficiently well-structured systems - this is our main result. As a simple illustration, we show how approximate controllability of a wave equation implies approximate controllability of a heat equation.
\end{abstract}

\section{Introduction}

Controllability is a fundamental property of dynamical systems and plays an important role in the analysis and design of control strategies. For infinite-dimensional systems, approximate controllability describes the possibility of steering the system from the origin to an arbitrarily close approximation of any prescribed target state in finite time \cite{CurtainZwart2020}. Understanding how this property depends on the underlying dynamics is a central question in control theory. In particular, it is important to determine whether controllability is preserved when a system is modified by interconnection with another system or by introducing dissipative mechanisms.
\smallskip

Port-Hamiltonian systems provide a natural framework for investigating such questions. Based on the underlying energy balance and the interconnection structure of a system, the port-Hamiltonian approach offers a systematic way of describing energy storage, energy transfer, and energy dissipation. It provides a unified mathematical framework for a broad class of finite- and infinite-dimensional systems, including systems governed by partial differential equations; see, for example, \cite{LeGorrec2005,vanDerSchaft2020,vanDerSchaftMaschke2002}.

A central concept in the construction of such systems is that of a closure relation. A closure relation specifies additional dissipation that determines the dynamics of the resulting system. In the context of port-Hamiltonian systems, Zwart, Le Gorrec, and Maschke \cite{ZwaH2016a} investigated how systems can be constructed from simpler hyperbolic systems using such interconnections. This approach was subsequently generalized by Schwenninger and Zwart \cite{SchF2014a}, who developed an operator-theoretic framework for generators with a closure relation on Banach spaces. More recently, Glück, Jacob, Meyer, Wyss, and Zwart \cite{GluJ2024a} studied stability via closure relations for port-Hamiltonian systems. 
\smallskip

In this article we investigate how the introduction of dissipation affects controllability. Our first main result shows that approximate controllability of the original system is, in general, insufficient to guarantee approximate controllability of the resulting target system. Thus, the introduction of a closure relation may alter the controllability properties in a way that cannot be inferred from the original system alone.
\smallskip

Our second main result establishes that the desired implication does hold under additional assumptions on the system operators. These assumptions provide sufficient conditions for transferring approximate controllability from the original system to the target system. \smallskip

This paper is organized as follows: The mathematical setup is presented in Section \ref{sec:Setting}. 
In Section \ref{sec:Conj} we present our initial conjecture about approximate controllability, and we show that it does not hold in general in Section \ref{sec:Counterexmpl} by means of a finite-dimensional counterexample. However, by considering additional assumptions, we prove that the conjecture holds in Section \ref{sec:A11=0}. This result is illustrated on an example in Section \ref{sec:example}. A conclusion and some perspectives can be found in Section \ref{sec:ConclPersp}. The appendix contains the proofs of several auxiliary results related to our main theorem. 

\section{Setting}\label{sec:Setting}

We consider the system
\begin{align}\label{eqn:system_ext}
\dot x(t) &= A_\ext x(t) +B_\ext u(t), \qquad t\ge 0,\\
y(t) &= B^*_\ext x(t)\nonumber
\end{align}
where the operators $A_\ext:\mathcal{D} (A_\ext) \subset X_1\times X_2 \to X_1\times X_2 $ and $B_\ext\in\mathcal{L}(U,X_1\times X_2 )$ are defined by \begin{align*}
    A_\ext:=\left[\begin{smallmatrix}
        A_{11} & A_{12} \\[1mm] A_{21} & 0
    \end{smallmatrix}\right], \quad B_\ext:=\left[\begin{smallmatrix}
        B \\[1mm] 0
    \end{smallmatrix}\right]
\end{align*} with domain
\begin{align*}
    \quad\mathcal{D} (A_\ext):=\left(\mathcal{D}(A_{11})\cap\mathcal{D}(A_{21})\right)\times \mathcal{D}(A_{21}).
\end{align*}
Here $X_1,X_2$ and $U$ are complex nonzero Hilbert spaces and the operators 
 $A_{11}:\mathcal{D}(A_{11})\subset X_1\to X_1$, $A_{12}:\mathcal{D}(A_{12})\subset X_2\to X_1$, $A_{21}:\mathcal{D}(A_{21})\subset X_1\to X_2$ and $B: U\to X_1$ are linear operators and $B\in\mathcal{L}(U,X_1)$. In what follows, we assume that $A_\ext$ generates a contraction strongly continuous ($C_0-$) semigroup. Thus, system \eqref{eqn:system_ext} possesses for every initial condition $x(0)\in X_1\times X_2$ and every input function $u\in \mathrm{L}^2((0,\infty);U)$ a unique mild solution $x(t)= \begin{pmatrix}x_1(t) & x_2(t)\end{pmatrix}^\top$, see \cite{CurtainZwart2020}. Further, system \eqref{eqn:system_ext} is a port-Hamiltonian system with Hamiltonian $H(x)=\|x\|^2_{X_1\times X_2}$.\smallskip

Moreover, let $R: X_2 \to X_2$ be a bounded linear operator which is positive definite, that is, $R = R^*$ and there exists $c > 0$ such that $\langle R x_2,x_2\rangle_{X_2}\geq c\Vert x_2\Vert_{X_2}^2$ for every $x_2\in X_2$. We now add to the system \eqref{eqn:system_ext}
a so-called \emph{closure relation} with $R$
and obtain the following port-Hamiltonian system with additional dissipation
\begin{align}\label{eqn:system_R}
    \dot x_1(t)&=A_Rx_1(t)+B_R
u(t), \qquad t\ge 0,\\
y(t) &=B^*_Rx_1(t),  \qquad t\ge 0,\nonumber
\end{align}
where 
\begin{align}
    A_R &= A_{11}+A_{12}RA_{21}, \quad B_R=B, \label{eq:A_R}\\
    \mathcal{D}(A_R)&=\left\{ x_1\in \mathcal{D}(A_{11})\cap\mathcal{D}(A_{21}) : RA_{21}x_1\in \mathcal{D}(A_{12}) \right\}. \nonumber
\end{align} 
 
Closure relations provide a systematic way to introduce
dissipation into a system by connecting a 
power port to a resistive element; see
\cite{GluJ2024a,Reis2021,
SchF2014a,vanDerSchaft2006}.
To explain this interpretation, we introduce the
power-conjugate effort and flow variables
$(e_1,f_1)$, $(e_2,f_2)$ and $(e_I,f_I)$, where
\begin{align}
    f_1 &=A_{11}e_1+A_{12}e_2 +Be_I,\nonumber\\
    f_2&=A_{21}e_1,\label{eqn:effort_flow}\\ 
    f_I &=B^*e_1.\nonumber
\end{align}
For the extended system \eqref{eqn:system_ext}, the dynamics $x$ then has to satisfy \eqref{eqn:effort_flow}
with $\left(\begin{smallmatrix}
    e_1 \\ e_2
\end{smallmatrix}\right)= x(t)$, $\left(\begin{smallmatrix}
    f_1 \\ f_2
\end{smallmatrix}\right)= \dot{x}(t)$, $e_I=u(t)$ and $f_I=y(t)$ for every $t\ge 0$.
The closure relation
\[
    e_2=Rf_2
\]
connects the  effort and flow variables $(e_2,f_2)$,
thereby eliminating this port. Moreover, the dynamics $x_1$ of system \eqref{eqn:system_R} has to satisfy \eqref{eqn:effort_flow} with $e_1=x_1(t)$, $f_1 =\dot{x}_1(t)$, $e_I=u(t)$ and $f_I=y(t)$ for every $t\ge 0$.
Since $R$ is positive definite,
the power absorbed through this port satisfies
\[
    \operatorname{Re}\langle f_2,e_2\rangle_{X_2}
    =\langle Rf_2,f_2\rangle_{X_2}
    \geq c\|f_2\|_{X_2}^2.
\]
If the extended operator $A_{\mathrm{ext}}$ is
dissipative, then for $e_I = 0$
\[
    \operatorname{Re}\langle f_1,e_1\rangle_{X_1}
    +
    \operatorname{Re}\langle f_2,e_2\rangle_{X_2}
    \leq 0.
\]
Consequently, we obtain for $e_I=0$
\begin{align*}
    \operatorname{Re}\langle A_Rx_1,x_1\rangle_{X_1}
    &=  \operatorname{Re}\langle(A_{11}+A_{12}RA_{21})e_1, e_1\rangle_{X_1}  \\
    &=  \operatorname{Re}\langle A_{11}e_1+A_{12}e_2, e_1\rangle_{X_1}  \\
     &=  \operatorname{Re}\langle f_1, e_1\rangle_{X_1}  \\
    &\leq
    -\operatorname{Re}\langle f_2,e_2\rangle_{X_2} \\
    &=-\langle f_2,R f_2\rangle_{X_2}\\
    &\leq 0.
\end{align*}
Thus, the closure relation introduces internal
dissipation while leaving the external input-output
pair unchanged. 

In \cite{GluJ2024a,PasH2026ax,SchF2014a, ZwaH2016a} the operator $A_\ext$ is given by $\left[
\begin{smallmatrix}
  A_1 \\[1mm]
  \begin{smallmatrix}
    A_{21} & 0
  \end{smallmatrix}
\end{smallmatrix}
\right]$ with domain $\left\{ \left(\begin{smallmatrix}
        x_1 \\x_2 
    \end{smallmatrix}\right)\in \mathcal{D}(A_1) : x_1\in \mathcal{D}(A_{21}) \right\}$ where $A_1:\mathcal{D}(A_1)\subset X_1\times X_2\to X_1$. In this paper, we assume that $A_1$ splits as $\begin{bmatrix}
        A_{11} &A_{12}
    \end{bmatrix}$.

Thanks to the boundedness of $B_R$, well-posedness of $\dot x_1(t)=A_Rx_1(t)+B_R
u(t), x_1(0) = x_{1_0}$, is equivalent to the question whether $A_R$ generates a $C_0$-semigroup. According to \cite[Thm.~2.2]{ZwaH2016a}, generation of a contraction semigroup by $A_\ext$ is sufficient to ensure that $A_R$ also generates a contraction semigroup.
\begin{lemma}\label{lem:Well-posedness}
    If $A_\ext$ generates a contraction semigroup on $X_1\times X_2$, then $A_R$ generates a contraction semigroup on $X_1$.
\end{lemma}

In the case where $A_{11}$ is everywhere defined and bounded, we have the following characterization for $A_\ext$ to be the generator of a contraction $C_0$-semigroup, see \cite[Thm.~2.8]{GluJ2024a}.
\begin{lemma}\label{lem:Adj}
    Assume that $A_{11}\in\mathcal{L}(X_1)$ is everywhere defined. Then, $A_\ext$ generates a contraction $C_0$-semigroup if and only if $A_{12}$ and $A_{21}$ are densely defined and closed with $A_{12}=-A_{21}^\ast$ and $A_{11}$ is dissipative. 
\end{lemma}

Before diving further in our main results, let us recall what we mean by approximate controllability.
\begin{definition}
    Let $X$ and $U$ be Hilbert spaces and consider an (unbounded) operator $A:\mathcal{D}(A)\subset X\to X$. Assume that $A$ generates a $C_0$-semigroup and that $B\in\mathcal{L}(U,X)$. The pair $(A,B)$ is said to be \emph{approximately controllable} if for every $\varepsilon>0$ and all $x_1\in X$ there exists a time $t>0$ and an input $u\in \mathrm{L}^2([0,t];U) $ such that $\left\| x_1-x(t)\right\|<\varepsilon$ where $x(t)$ is the mild solution of $\dot x(t) =Ax(t)+Bu(t), x(0)=0$ at time $t$.
\end{definition}

\section{The conjecture}\label{sec:Conj}
Following the same approach as for the generation of semigroup and the stability, our tentative idea was to think that approximate controllability transferred from $(A_\ext, B_\ext)$ to $(A_R,B_R)$. This is summarized in the following conjecture.
\begin{conjecture}\label{conj}
 Assume that the pair $(A_\ext, B_\ext)$ is approximately controllable. Then, the pair $(A_R,B_R)$ obtained from $(A_\ext,B_\ext)$ and the closure relation induced by $R$ is approximately controllable.   
\end{conjecture}

As we will see later in this section, Conjecture \ref{conj} is not true in general. From \cite{GluJ2024a}, it is known that \begin{align*}
    (\lambda-A_R)^{-1} Bu = \left[\begin{smallmatrix}
        I & 0 
    \end{smallmatrix}\right] (\lambda-(A_\ext+P_\lambda))^{-1} \left[\begin{smallmatrix}
        Bu \\[1mm] 0
    \end{smallmatrix}\right]
\end{align*} for $\lambda\in\rho(A_R)\cap \rho(A_\ext)$ and $u\in U$, where $P_\lambda=\left[\begin{smallmatrix}
    0 & 0 \\[1mm] 0& -R^{-1} + \lambda
\end{smallmatrix}\right]$. Our idea to prove Conjecture \ref{conj} was based on the characterization of approximate controllability by the resolvent operator, see e.g.~\cite[Lem.~6.2.24]{CurtainZwart2020} 
\begin{lemma}\label{lem:ControlViaResolvent}
    Let $X, U$ be Hilbert spaces and let $A$ be the generator of a contraction $C_0$-semigroup on $X$ and $B\in\mathcal{L}(U,X)$. The pair $(A,B)$ is approximately controllable if and only if 
    \begin{equation*}
    \overline{\spn_{u\in U, \lambda>0} (\lambda-A)^{-1}Bu} =X.
    \end{equation*}
\end{lemma}

\section{A counterexample for the general case}\label{sec:Counterexmpl}
The following simple counterexample shows that Conjecture \ref{conj} does not hold in general.
	
    \begin{example}\label{ex:main}
    We consider $X_1 = \mathbb{R}^2, X_2 = \mathbb{R}$ and $U = \mathbb{R}$, and we set
	\begin{align*}
		A_\ext 
		:= 
		\left(
			\begin{smallmatrix}
				0 & -1 & 1 \\[1mm]
				1 &  0 & 1 \\[1mm]
				-1 & -1 & 0
			\end{smallmatrix}
		\right), \quad B_\ext:=\left(
			\begin{smallmatrix}
				1\\[1mm]
                0\\[1mm]
				0
			\end{smallmatrix}
		\right), \quad R := I.
	\end{align*}
	Lemma \ref{lem:Adj} implies that $A_\ext$ generates a contraction semigroup. The following Kalman rank type condition shows that $(A_\ext, B_\ext)$ is controllable\footnote{For finite-dimensional systems, approximate and exact controllability are equivalent concepts. Therefore we say that $(A_\ext, B_\ext)$ is controllable in the finite-dimensional setting.}, see e.g.~\cite[Chap.~6]{HinD2026b2}:
    \begin{align*}
        \rank\begin{bmatrix}
            B_\ext & A_\ext B_\ext & A_\ext^2B_\ext
        \end{bmatrix}=\rank\left(\begin{smallmatrix}
            1 & 0 &  -2\\[1mm]
            0 & 1 &  -1\\[1mm]
            0 & -1 & -1
        \end{smallmatrix}\right)=3.
    \end{align*}
    Moreover, we have that
	\begin{align*}
		A_R 
		= 
		\left(\begin{smallmatrix}
			-1 & -2 \\[1mm] 
			 0 & -1
		\end{smallmatrix}\right), \quad B_R=\left(\begin{smallmatrix}
		    1\\[1mm] 0
		\end{smallmatrix}\right)
		.
	\end{align*}
    Applying the Kalman rank condition to $(A_R,B_R)$, we get that \begin{align*}
        \rank \begin{bmatrix}
            B_R & A_R B_R
        \end{bmatrix}=\rank \left(\begin{smallmatrix}
            1 & -1 \\[1mm]
            0 &0
        \end{smallmatrix}\right)=1,
    \end{align*}    
   which shows that $(A_R,B_R)$ is not controllable.
    \end{example}

\section{The case $A_{11}=0$}\label{sec:A11=0}
In this section we consider the case where $A_{11}=0$, which leads to $A_R=A_{12}RA_{21}$. Without adapting our current setting, we start by showing that Conjecture \ref{conj} is still not valid.

\begin{example}\label{ex:CountExampl}
Let $X_1 = X_2 = \mathbb{R}^2$ and $U = \mathbb{R}$. We define
\begin{align*}
A_\ext &:=
\left(\begin{smallmatrix}
0 & 0 & 1 & 1\\[1mm]
0 & 0 & 0 & 1\\[1mm]
-1 & 0 & 0 & 0\\[1mm]
-1 & -1 & 0 & 0
\end{smallmatrix}\right), \quad B_\ext := \left(\begin{smallmatrix}
    1 \\[1mm] 0\\[1mm] 
    0\\[1mm] 0
\end{smallmatrix}\right),\quad R :=
\left(\begin{smallmatrix}
2 & -1\\[1mm]
-1 & 1
\end{smallmatrix}\right).    
\end{align*}
As it has been done in Example \ref{ex:main}, a Kalman rank condition can be used to show that $(A_\ext, B_\ext)$ is controllable, but that $(A_R,B_R)$ is not.
\end{example}

Example \ref{ex:CountExampl} shows that for $A_{11}=0$ approximate controllability does not transfer from $(A_\ext,B_\ext)$ to $(A_R,B_R)$, which will force us adapting the setting. We now consider the operators
\begin{equation}
    \tilde{A}_\ext=\left[\begin{smallmatrix}
        0 & A_{12} \\[1mm] A_{21} &0
    \end{smallmatrix}\right]\left[\begin{smallmatrix}
        I & 0 \\[1mm] 0 & R
    \end{smallmatrix}\right], \quad \tilde{B}_\ext=\left[\begin{smallmatrix}
        B \\[1mm] 0
    \end{smallmatrix}\right]
    \label{eq:Aext_Bext_t}
\end{equation}
where 
\begin{equation}
\mathcal{D}(\tilde{A}_\ext)=\mathcal{D}(A_{21})\times \left\{x_2\in X_2 : Rx_2\in \mathcal{D}(A_{12}) \right\}
\label{eq:Dom_Aext_t}
\end{equation}
in which $R$ is a positive definite operator on $X_2$ and $A_{12}$ and $A_{21}$ are closed and densely defined linear operators, see Section \ref{sec:Setting}. According to Lemma \ref{lem:Adj}, $\left[\begin{smallmatrix}
        0 & A_{12} \\[1mm] A_{21} &0
    \end{smallmatrix}\right]$ generates a contraction $C_0$-semigroup if and only if $A_{12}^\ast=-A_{21}$. We make this assumption from now on. Thus, $\tilde A_\ext$ generates a contraction semigroup in the so-called energy norm, induced by the inner product 
    \begin{equation*} \left\langle x,\tilde{x}\right\rangle_R := \left\langle x,\left[\begin{smallmatrix}
        I & 0 \\[1mm] 0 & R
    \end{smallmatrix}\right]\tilde x\right\rangle_{X_1\times X_2} \quad x,\tilde x\in X_1\times X_2.
    \end{equation*}
    Observe also that the system $\dot x(t) =\tilde A_\ext x(t) +\tilde B_\ext u(t)$ with the closure relation induced by $I_{X_2}$ leads to $\dot x_1(t)=A_Rx_1(t)+B_Ru(t)$ where $A_R = A_{12}RA_{21} = -A_{12}RA_{12}^*$ and $B_R = B$. The contraction semigroup generation of $A_R$ is a consequence of Lemma \ref{lem:Well-posedness}.

Before stating our main result, we need to recall the following characterization of approximate controllability, due to Fattorini, see \cite{Fattorini_Contr} and \cite[Thm.~1.3]{Badra_Takahashi}.

\begin{lemma}\label{lem:Fattorini}
    Let $X,U$ be Hilbert spaces and let $B \in\mathcal{L}(U,X)$. 
    Let $A$ be the generator of a $C_0$-semigroup on $X$ whose spectrum consists of isolated eigenvalues with finite algebraic multiplicity and assume that the corresponding eigenvectors are complete in $X$. 
    Then $(A,B)$ is approximately controllable if and only if
    \begin{equation*}
        \ker(\lambda I-A^*)\, \cap\, \ker(B^*) = \{0\}
    \end{equation*}
    holds for all $x \in \mathcal{D}(A^*)$ and all $\lambda\in\mathbb{C}$.
\end{lemma}

Before presenting our main result, we recall the following: a point $\lambda\in\sigma_p(A)$ is an \emph{isolated eigenvalue} of $A$ if there exists $r_\lambda > 0$ such that $\sigma(A)\cap\{\mu\in\mathbb{C}, \vert\mu-\lambda\vert < r_\lambda\} = \{\lambda\}$ which means that there exists a disk around $\lambda$ which contains no other spectral point.

Our main positive result concerning approximate controllability is based on Lemma~\ref{lem:Fattorini} and is summarized in the following theorem.

\begin{theorem}\label{thm:A11=0}
  Consider the operators $\tilde{A}_\ext, \tilde{B}_\ext$ from \eqref{eq:Aext_Bext_t}--\eqref{eq:Dom_Aext_t} with $A_{12}^*=-A_{21}$. Then, approximate controllability of $(\tilde{A}_{\ext},\tilde{B}_\ext)$ implies approximate controllability of $(A_R,B_R)$. If in addition $\ker A_{12} = \{0\}$ and $A_R$ has compact resolvent\footnote{Equivalent conditions for $A_R$ to have compact resolvent are established in Proposition~\ref{prop:CompRes}.}, then approximate controllability of $(A_R,B_R)$ implies approximate controllability of $(\tilde{A}_\ext,\tilde{B}_\ext)$.
\end{theorem}

\begin{proof}
For simplicity in the notations, let us introduce $J := A_{12}$. With this definition, $\tilde{A}_{\ext} = \left[\begin{smallmatrix}0 & JR\\[1mm] -J^* & 0\end{smallmatrix}\right]$ and $A_R = -JRJ^*$.  
We start by showing that approximate controllability of $(\tilde{A}_\ext,\tilde{B}_\ext)$ implies approximate controllability of $(A_R,B_R)$, for which we shall apply Lemma \ref{lem:ControlViaResolvent}. For $\lambda>0$, observe that the resolvent operator of $\tilde A_\ext$ is given by
\begin{align*}
       & (\lambda-\tilde A_\ext)^{-1}=\left[\begin{smallmatrix}
\lambda(\lambda^2+JRJ^*)^{-1}
&
JR(\lambda^2+J^*JR)^{-1}
\\[1mm]
-J^*(\lambda^2+JRJ^*)^{-1}
&
\lambda(\lambda^2+J^*JR)^{-1}
\end{smallmatrix}\right].
    \end{align*}
Consequently, we have that \begin{align*}
    \left[\begin{smallmatrix}
        I & 0
    \end{smallmatrix}\right](\lambda-\tilde A_\ext)^{-1} \left[\begin{smallmatrix}
        B \\[1mm] 0
    \end{smallmatrix}\right]= \lambda (\lambda^2-A_R)^{-1}B.
\end{align*}
Since $(\tilde A_\ext,\tilde B_\ext)$ is approximately controllable, we have \begin{align*}
     X_1 
    &= \left[\begin{smallmatrix}
        I & 0
    \end{smallmatrix}\right] \overline{\spn_{u\in U, \lambda>0} (\lambda-\tilde A_\ext)^{-1}\tilde B_\ext u}\\
    &=\overline{\spn_{u\in U, \lambda>0} \lambda (\lambda^2-A_R)^{-1}B_R u}\\
    &=\overline{\spn_{u\in U, \lambda>0} (\lambda-A_R)^{-1}B_R u}.
\end{align*}
According to Lemma \ref{lem:ControlViaResolvent}, the previous calculation shows that $(A_R,B_R)$ is approximately controllable.

Assume now that $(A_R,B_R)$ is approximately controllable. Because $A_R = -JRJ^*$ is self-adjoint, nonpositive and has compact resolvent, \cite[Lem.~3.2.9]{CurtainZwart2020} implies  
that the eigenfunctions of $A_R$ form an orthonormal basis in $X_1$. Moreover, compact resolvent of $A_R$ implies that $\sigma(A_R)$ is composed of isolated eigenvalues only. We shall now apply Lemma~\ref{lem:Fattorini}. 
Let $x=\left(\begin{smallmatrix}
    x_1 \\[1mm] x_2
\end{smallmatrix}\right)\in\mathcal{D}(\tilde{A}_\ext)$ and $\lambda\in\mathbb C$ be such that $\tilde A_\ext^\ast x=\lambda x$ and $\tilde B_\ext^\ast x =0$.
From the structure of $\tilde A_\ext$ we have that $-J x_2=\lambda x_1$ and $RJ^\ast x_1=\lambda x_2$. 
Hence, we have that $-JRJ^*x_1=-\lambda Jx_2$. With the definition of $A_R$ we conclude that $A_R x_1=\lambda^2 x_1$. Furthermore $\tilde B_\ext^\ast x =0$ implies $B_R^\ast x_1=0$. Thus Lemma \ref{lem:Fattorini} implies $x_1=0$. Going back to equation $-Jx_2 = \lambda x_1$ with $x_1=0$ and the assumption $\ker J = \{0\}$ shows that $x_2=0$. We conclude that $x=0$. In order to apply Lemma~\ref{lem:Fattorini} to $(\tilde{A}_\ext,\tilde{B}_\ext)$, it remains to show that 
 $\tilde{A}_\ext$ has compact resolvent and   
that the corresponding eigenfunctions are complete in $X_1\times X_2$ to conclude the proof. 
Due to the assumptions that $A_R$ has compact resolvent and that $\ker A_{12}$ is finite-dimensional, $\tilde A_\ext $ has a compact resolvent, see Proposition~\ref{prop:CompRes}. Noting that $\tilde{A}_\ext^* = -\tilde{A}_\ext$ on $X_1\times X_2$ equipped with the energy inner product $\langle\cdot,\cdot\rangle_R$,  
\cite[Lem.~3.2.9]{CurtainZwart2020} implies that the corresponding eigenfunctions are complete in $X_1\times X_2$.
\end{proof}

\begin{remark}\label{rem:Main}
    \begin{itemize}
        \item In the case  $R=I$,  $(\tilde A_\ext,\tilde B_\ext)$ equals  $(A_\ext,B_\ext)$. 
        Hence, Theorem \ref{thm:A11=0} shows that Conjecture \ref{conj} is true when $A_{11} = 0$, $A_{12}^\ast=-A_{21}$ and $R = I$.
        \item Note that $A_\ext$, $R$ and $B_\ext$ from Example \ref{ex:CountExampl} do not fit into the setting of Theorem \ref{thm:A11=0}. Here we have 
        \begin{equation*}
        \tilde A_\ext= \left(\begin{smallmatrix}
            0 & 0 & 1 & 0 \\[1mm]
            0 & 0 & -1 & 1 \\[1mm]
            -1 & 0&0&0 \\[1mm]-1&-1&0&0
        \end{smallmatrix}\right), \quad \tilde B_\ext=\left(\begin{smallmatrix}
            1 \\[1mm]0\\[1mm]0\\[1mm]0
        \end{smallmatrix}\right).
        \end{equation*}
    Thus, the pair $(\tilde A_\ext, \tilde B_\ext$) is not controllable.
    \end{itemize}
\end{remark}

It is worth noting that without the assumption that $\ker A_{12} = \{0\}$, the converse of Theorem \ref{thm:A11=0} might not be true, which is highlighted in the following example. 

\begin{example}\label{ex:Converse}
    Consider $\tilde A_\ext$ with $A_{12}=\left(\begin{smallmatrix}-1 & 0\end{smallmatrix}\right)$, 
$R=\left(\begin{smallmatrix}1&0\\[1mm] 0&1\end{smallmatrix}\right)$, 
$A_{21}=\left(\begin{smallmatrix}1\\[1mm] 0\end{smallmatrix}\right)$ as well as $\tilde B_\ext = \left(\begin{smallmatrix}1\\[1mm] 0\end{smallmatrix}\right)$ with
$B=1$. It is clear that $\ker A_{12}$ contains any vector of the form $\left(\begin{smallmatrix}0\\[1mm] a\end{smallmatrix}\right)$, $a\in\mathbb{C}$, so that $\ker A_{12} \neq \{0\}$. In that case, $A_R=-1$ and $B_R=1$, which implies that $(A_R,B_R)$ is controllable, but it is easy to see that $(\tilde A_\ext, \tilde B_\ext)$ is not controllable.
\end{example}

\section{Example}\label{sec:example}
In this section we apply Theorem \ref{thm:A11=0} on a wave to heat example. More precisely, we show that the equivalent conditions characterizing the approximate controllability of a vibrating string with distributed control are the same conditions which characterize the approximate controllability of a heat equation obtained by considering the vibrating string with a closure relation. We start by considering a vibrating string whose vertical displacement, denoted by $w(\xi,t)$, $\xi\in (0,1), t\geq 0$, is subject to the following controlled one-dimensional wave equation \begin{equation}
    \frac{\partial^2 w}{\partial t^2}(\xi,t)= \frac{\partial^2 w}{\partial \xi^2}(\xi,t) + b(\xi)u(t).
    \label{eq:wave}
\end{equation}
We further assume that the string is clamped at its left boundary and that no force is applied to the right boundary, resulting in the following boundary conditions 
\begin{equation*}
    \frac{\partial w}{\partial t}(0,t)=0, \qquad \frac{\partial w}{\partial \xi}(1,t)=0, \qquad t\geq 0.
\end{equation*}    
In \eqref{eq:wave}, we assume that the scalar input is distributed in the domain thanks to the function $b$. In particular, we assume that $b\in\mathrm L^2(0,1)$ and $u(t)\in\mathbb C$. By defining the energy variables $x_1=\frac{\partial w}{\partial t} $ and $x_2=\frac{\partial w}{\partial \xi}$, \eqref{eq:wave} is equivalent to 
\begin{equation}
    \frac{\partial x}{\partial t}(\xi,t)= \begin{pmatrix}
        0 & \frac{\partial}{\partial\xi} \\ \frac{\partial}{\partial\xi} & 0
    \end{pmatrix}x(\xi,t)+ \begin{pmatrix}
        b(\xi) \\0
    \end{pmatrix}u(t),
    \label{eq:Wave_pHs}
\end{equation}
with boundary conditions $x_1(0,t)=0$ and $x_2(1,t)=0$, $t\geq 0$. It is easy to see that \eqref{eq:Wave_pHs} can be written in an abstract way thanks to the operators $\tilde{A}_\ext$ and $\tilde{B}_\ext$ from \eqref{eq:Aext_Bext_t}, where
\begin{align*}
    A_{12}f &= \frac{\mathrm{d}f}{\mathrm{d}\xi},\quad     \mathcal{D}(A_{12}) = \left\{ f\in \mathrm{H}^1(0,1) : f(1)=0 \right\},
\end{align*}
$R = I_{X_2}$ and $B:U\to X_1, Bu = bu$, where $U = \mathbb{C}$ and $X_1 = X_2 = \mathrm{L}^2(0,1)$. Further, the operator $A_{21}$ is given by
\begin{align*}
    A_{21}f &= \frac{\mathrm{d}f}{\mathrm{d}\xi},\quad \mathcal{D}(A_{21}) = \left\{ f\in \mathrm{H}^1(0,1) : f(0)=0 \right\}
\end{align*}
which implies that $A_{21} = -A_{12}^\ast$. Hence, the operator $\tilde A_\ext$ generates a contraction $C_0$-semigroup, see Lemma \ref{lem:Adj}.

Moreover, it can be shown that $\tilde{A}_\ext$ is a Riesz-spectral operator in the sense of \cite[Def.~3.2.6]{CurtainZwart2020}. The eigenvalues of $\tilde{A}_\ext$, denoted by $\{\lambda_n\}_{n\in\mathbb{Z}}$, and the corresponding Riesz basis of eigenfunctions, denoted by $\{\phi_n\}_{n\in\mathbb{Z}}$, are given by 
\begin{align*}
\lambda_n = \mathrm i \kappa_n,\quad \phi_n(\xi)&=\left(\begin{smallmatrix}
    \sin(\kappa_n \xi)\\[1mm] \mathrm i \cos(\kappa_n \xi)
\end{smallmatrix}\right), \quad \xi\in[0,1],
\end{align*}
respectively, with $\kappa_n=\left( n+ \frac{1}{2} \right) \pi$. The operator $\tilde A_\ext$ is skew-adjoint, which implies that the eigenvalues and the eigenfunctions of $\tilde A_\ext^\ast$ are given by $\overline{\lambda_n}=-\mathrm i \kappa_n$ and $\psi_n = \phi_n, n\in\mathbb{Z}$, respectively.

To characterize the controllability of $(\tilde{A}_\ext,\tilde{B}_\ext)$, we apply \cite[Thm.~6.3.6]{CurtainZwart2020}. It follows that $(\tilde A_\ext, \tilde B_\ext)$ is approximately controllable if and only if $\left\langle \left(\begin{smallmatrix}
    b \\[1mm] 0
\end{smallmatrix}\right), \psi_n \right\rangle_{\mathrm L^2(0,1) \times\mathrm L^2(0,1) } \neq 0$ for every $n\in\mathbb Z$.
Because of $\kappa_{-n-1}=-\kappa_n$, $n\in\mathbb{N}$, the previous condition is equivalent to \begin{equation}
    \left\langle b(\cdot), \sin(\kappa_n\cdot) \right\rangle_{\mathrm L^2(0,1)} \neq 0, \quad n\in \mathbb N_0.
    \label{eq:ContrTest_Wave}
\end{equation}

Considering \eqref{eq:Wave_pHs} together with the closure relation induced by $R = I_{X_2}$ yields the following controlled heat equation
\begin{align*}
     \frac{\partial x_1}{\partial t}(\xi,t) = \frac{\partial^2 x_1}{\partial\xi^2}(\xi,t)+b(\xi) u(t)
\end{align*} 
with boundary conditions $x_1(0,t)=0$ and $\frac{\partial x_1}{\partial\xi}(1,t)=0$, $\xi\in[0,1]$, $t\geq 0$. According to \eqref{eq:A_R}, the operators $A_R$ and $B_R$ are given by
\begin{align*}
    A_R f &= \frac{\mathrm{d}^2f}{\mathrm{d}\xi^2},\, D(A_R) = \{f\in\mathrm{H}^2(0,1), f(0) = 0 = \frac{\mathrm{d}f}{\mathrm{d}\xi}(1)\}\\
    B_Ru &= bu, 
\end{align*}
respectively. It is easy to verify that $A_R$ is a self-adjoint Riesz-spectral operator with an orthonormal basis of eigenvectors. The eigenvalues of $A_R$ are given by $\tilde{\lambda}_n = -\kappa_n^2, n\in\mathbb{N}_0$ and the corresponding eigenfunctions are expressed as $\tilde{\phi}_n(\xi) = \sqrt{2}\sin(\kappa_n\xi), \xi\in [0,1], n\in\mathbb{N}_0$, where $\kappa_n = \left( n+\frac{1}{2}\right)\pi$. Applying the same controllability test as the one used for the wave equation, that is, \cite[Thm.~6.3.6]{CurtainZwart2020}, yields that $(A_R,B_R)$ is approximately controllable if and only if $\langle b(\cdot),\sin(\kappa_n\cdot)\rangle_{\mathrm{L}^2(0,1)}\neq 0$ for every $n\in\mathbb{N}_0$. This characterization is exactly the same as \eqref{eq:ContrTest_Wave}, which already means that $(\tilde{A}_\ext,\tilde{B}_\ext)$ is approximately controllable if and only if $(A_R,B_R)$ is approximately controllable. We shall now confirm this statement with Theorem \ref{thm:A11=0}. Thanks to our setting, the last assumption that is needed to apply Theorem \ref{thm:A11=0} is $\ker A_{12} = \{0\}$. Take $f\in\ker A_{12}$. Hence, we have $\frac{\mathrm{d}f}{\mathrm{d}\xi} = 0$, which yields $f(\xi) = c$ on $[0,1]$, where $c$ is a constant. Thanks to $f(1) = 0$, we get that $c = 0$, yielding $\ker A_{12} = \{0\}$. Theorem \ref{thm:A11=0} then corroborates that approximate controllability of $(\tilde{A}_\ext,\tilde{B}_\ext)$ is equivalent to approximate controllability of $(A_R,B_R)$.

\section{Conclusion and perspectives}\label{sec:ConclPersp}
By using different counterexamples, we have shown that controllability of a linear system obtained by means of a closure relation cannot be guaranteed by the controllability of the original system. Nevertheless, we have identified a setting in which controllability of the nominal and the target systems are equivalent, which is the main result of this paper. We applied this result to show that approximate controllability of a wave equation is equivalent to approximate controllability of a heat equation, where a closure relation together with the port-Hamiltonian representation of the wave have been used.

Further research could aim at investigating other notions of controllability in this context. Extending the proposed approach to boundary control could also be considered as perspective. 

The approach based on closure relations could also be investigated in the context of linear quadratic optimal control or feedback stabilization.

\appendix[Compact resolvents]\label{sec:appendix}
\setcounter{theorem}{0}
\renewcommand{\thetheorem}{A.\arabic{theorem}}

In the reverse direction of Theorem~\ref{thm:A11=0} we assume that $A_R$ has compact resolvent and $\ker A_{12} = \{0\}$. 
The slightly more general property that $A_R$ has compact resolvent and $\dim \ker A_{12} < \infty$, can be characterized in various different ways:

\begin{proposition}\label{prop:CompRes}
    The following are equivalent: \begin{enumerate}
        \item $A_\ext$ has compact resolvent.
        \item There exists a positive definite $R\in \mathcal{L}(X_2)$ such that $\tilde A_\ext$ has compact resolvent.
        \item For all positive definite $R\in \mathcal{L}(X_2)$, $\tilde A_\ext$ has compact resolvent.
        \item There exists a positive definite $R\in \mathcal{L}(X_2)$ such that $ A_R$ has compact resolvent, and $\dim \ker A_{12} < \infty$.
        \item For all positive definite $R\in \mathcal{L}(X_2)$ the operator $A_R$ has compact resolvent, and $\dim \ker A_{12} < \infty$.
        \item $\mathcal D(A_{12})$ equipped with the graph norm embeds compactly into $X_2$ and $\dim \ker A_{21} < \infty$.
        \item $\mathcal D(A_{21})$ equipped with the graph norm embeds compactly into $X_1$ and $\dim \ker A_{12} < \infty$.
    \end{enumerate}
\end{proposition}

To prove Proposition \ref{prop:CompRes} we need the following three lemmas.

\begin{lemma}\label{lem:HelpWithCompRes}
    Consider a closed and densely defined linear operator $J:\mathcal{D}(J)\subset X_2 \to X_1$. 
    \begin{enumerate}
        \item 
        If $\mathcal{D}(J)$ equipped with the graph norm $\| \cdot\|_{\mathcal{D}(J)}^2:= \|\cdot\|^2_{X_2}+ \|J\cdot \|^2_{X_1}$ is compactly embedded in $X_2$, then $\ker J$ is finite dimensional.
        
        \item 
        Assume now that $X_1 = X_2$ and that $J$ has non-empty resolvent set. 
        Then $J$ has compact resolvent if and only if $\mathcal{D}(J)$ equipped with the graph norm is compactly embedded in $X$.
    \end{enumerate}
\end{lemma}

\begin{proof}
    For all $x\in\ker J$ one has $\|x\|_{\mathcal{D}(J)}=\|x\|_{X_2}$. 
    Since the unit ball of $\ker J$ is relatively compact in $X_2$, it follows that it is also relatively compact in $\ker J$.
    Thus, $\ker J$ is finite-dimensional. 
    The second assertion follows from e.g.~\cite[Prop.~4.25]{EngK2000b}. 
\end{proof}

\begin{lemma}\label{lem:CompRes}
    Let $J:\mathcal{D}(J)\subset X_2 \to X_1$ be a closed, densely defined linear operator. 
    The following are equivalent:
    \begin{enumerate}
        \item $\mathcal{D}(J^\ast)$ with the graph norm embeds compactly into $X_1$.
        \item  There exists a positive definite $R\in \mathcal{L}(X_2)$ such that $JRJ^*$ has compact resolvent.
        \item For all positive definite $R\in \mathcal{L}(X_2)$, $JRJ^*$ has compact resolvent.
    \end{enumerate}            
\end{lemma}

\begin{proof}
        To show that 
        {1)$\implies$3)}, let $R\in \mathcal{L}(X_2)$ be  positive definite.
        Firstly we show that the embedding $\tilde\iota_{\mathcal{D}(JRJ^*)}:\mathcal{D}(JRJ^*)\to \mathcal{D}(J^*) $ is continuous. As the embedding of $\mathcal{D}(J^\ast)$ into $X_1$ is compact and thus continuous, there exists $c>0$ such that for all $x\in \mathcal{D}(JRJ^*)$, we have 
        \begin{align*}
            \frac{\left\|x \right\|_{X_1} \left\| x  \right\|_{\mathcal{D}(J^\ast)}}{c}  
            &\le \left\| x  \right\|_{\mathcal{D}(J^\ast)}^2 \\
            &= \left\|R^{-\frac{1}{2}} R^\frac{1}{2}  J^\ast x \right\|_{X_2}^2+ \left\|x \right\|_{X_1}^2 \\
            &= \left\|R^{-\frac{1}{2}}  \right\|^2 \left| \left\langle  x, JR J^\ast x \right\rangle_{X_1}\right| + \left\|x \right\|_{X_1}^2 \\
            &\le  \left\|R^{-\frac{1}{2}}   \right\|^2 \left\|x \right\|_{X_1}  \left\|JR J^\ast x \right\|_{X_1} + \left\|x \right\|_{X_1}^2.
        \end{align*}
        Dividing by $\left\|x \right\|_{X_1}$ shows that $\tilde\iota_{\mathcal{D}(JRJ^*)}$ is bounded. Thus the embedding $\iota_{\mathcal{D}(JRJ^*)}:\mathcal{D}(JRJ^*)\to X_1$ is compact, because we have $\iota_{\mathcal{D}(JRJ^*)}= \iota_{\mathcal{D}(J^\ast)} \circ \tilde\iota_{\mathcal{D}(JRJ^*)} $. Hence, Lemma~\ref{lem:HelpWithCompRes} implies that $JRJ^\ast$ has compact resolvent.
    The implication {3)$\implies$2)} is clear. To show the implication {2)$\implies$1)}, consider the quadratic form associated to $JRJ^\ast$, given by $a(x,\tilde x):=\left\langle J^\ast x,  RJ^\ast \tilde x \right\rangle_{X_2}$, $x,\tilde x\in \mathcal{D}(J^\ast)$. Since $R$ is positive definite, there exists $c>0$ such that for all $x\in \mathcal D (J^\ast)$, we have $\mathrm{Re}\, a(x,x) + c\| x \|_{X_2} \ge c \|x\|_{\mathcal D(J^\ast)}$. Thus \cite[Ex.~6.8]{AreW2026pp} implies that the embedding $\iota_{\mathcal D(J^\ast)}$ is compact.
\end{proof}

\begin{lemma}\label{lem:CompResDual}
        Let $J:\mathcal{D}(J)\subset X_2 \to X_1$ be a closed and densely defined linear operator. 
        Then the following are equivalent: 
        \begin{enumerate}
            \item $\mathcal{D}(J)$ equipped with the graph norm embeds compactly in $X_2$ and $\ker J^\ast$ is finite dimensional.
            \item $\mathcal{D}(J^\ast)$ equipped with the graph norm embeds compactly in $X_1$ and $\ker J$ is finite dimensional.
        \end{enumerate}
\end{lemma}
\begin{proof}
    Assume that $\mathcal D (J)$ is compactly embedded in $X_2$ and that $\ker J^\ast$ is finite-dimensional. Lemma~\ref{lem:HelpWithCompRes} yields that $\ker J$ is finite-dimensional.
    It remains to show that $\mathcal{D}(J^\ast)$ is compact in $X_1$. 
    Let $(x_n)_n\subset\mathcal D(J^\ast)$ be a weak null sequence. 
    We need to show that $(x_n)_n$ converges to zero in norm. 
    Since norm continuous linear maps are also weak-to-weak continuous, we have $J^\ast x_n \rightharpoonup 0$ in $X_2$ and $x_n \rightharpoonup 0$ in $X_1$, where ``$\rightharpoonup$" stands for the weak convergence.
    Since adjoints of compact operators are compact by Schauder's Theorem, see e.g.~\cite[Thm~6.4]{Brezis2010}, $\iota_{\mathcal{D}(J)}^\ast:X_2\to \mathcal{D}(J)^\ast$ is also compact%
    \footnote{
        We identify $X_2$ with its own dual space; 
        however, in what follows, it is useful to distinguish between $\mathcal{D}(J)$ and its dual $\mathcal{D}(J)^*$.
    }.
    Thus, $J^\ast x_n \to 0$ in $\mathcal{D}(J)^\ast$ and for all bounded sequences $(y_n)_n\subset \mathcal{D}(J)$, we have
    \begin{align}
        \label{eq:test-against-bdd-sequence}
        \left\langle x_n, Jy_n\right\rangle_{X_1}= \left\langle J^\ast x_n, y_n\right\rangle_{X_2}= \iota_{\mathcal{D}(J)}^\ast(J^\ast x_n(y_n)) \to 0.
    \end{align}
    
    We will now show that $\ran J$ is closed and has finite codimension.%
    \footnote{
        As $\ker J$ is also finite-dimensional, $J:\mathcal{D}(J)\to X_1$ is a Fredholm operator.
    }
    We split $\mathcal D(J) = \ker J \oplus (\ker J)^\perp$, where the orthogonal complement is taken with respect to the inner product on $\mathcal{D}(J)$ that induces the graph norm of $J$.
    We show that $J\vert_{(\ker J)^\perp}$ is bounded from below; 
    this will imply that this operator has closed range and thus, so has $J$.
    Assume that $J\vert_{(\ker J)^\perp}$ is not bounded from below. 
    Then there exists a sequence $(z_n)_n$ in $(\ker J)^\perp$ with $\| z_n\|_{\mathcal{D}(J)}=1$ and $Jz_n\to 0$. 
    After replacing $(z_n)_n$ with a subsequence we may assume that it converges weakly in $\mathcal{D}(J)$, and the injectivity of $J$ on $(\ker J)^\perp$ implies that the weak limit is $0$. 
    Since $\iota_{\mathcal{D}(J)}$ is compact, $\iota_{\mathcal{D}(J)} z_n \to 0$ with respect to the norm in $X_2$. 
    Since $\mathcal{D}(J)$ is endowed with the graph norm, 
    this implies that $(z_n)_n$ is norm convergent to $0$ in $\mathcal{D}(J)$. 
    This contradicts $\|z_n\|_{\mathcal{D}(J)}=1$. 
    So the range of $J$ is  closed.
    
    Since $\ker J^\ast = (\ran J)^\perp$ is isomorphic to $X_2 /\overline{\ran J} =X_2 / {\ran J}$, the assumption $\dim \ker J^\ast < \infty$ implies that $\ran J$ has finite codimension.

    Let $P_{\ran J}$ be the orthogonal projection on $\ran J$. 
    Then
    \begin{align*}
        \left\langle x_n, x_n \right\rangle_{X_1}
        =  
        \left\langle x_n, (I-P_{\ran J}) x_n \right\rangle_{X_1}   + \left\langle x_n, P_{\ran J} x_n \right\rangle_{X_1}.     
    \end{align*}    
    The first summand converges to $0$ since $(x_n)_n$ converges weakly to $0$ and $(I-P_{\ran J})$ has finite rank. 
    To see that the second summand converges to $0$, we use that $J$ is an open map onto its range as its range is closed: 
    As $(P_{\ran J} x_n)$ is a bounded sequence on $\ran J$, the openness of $J$ implies that there exists a bounded sequence $(y_n)_n$ in $\mathcal{D}(J)$ such that $Jy_n =  P_{\ran J}x_n$ for each $n$. 
    Hence, \eqref{eq:test-against-bdd-sequence} yields
    \begin{align*}
        \left\langle x_n, P_{\ran J} x_n \right\rangle_{X_1} 
        = 
        \left\langle x_n, Jy_n \right\rangle_{X_1} 
        \to 
        0.
        \tag*{\QED}
    \end{align*}
    \let\QED\relax
\end{proof}
Now we are able to prove Proposition~\ref{prop:CompRes}:

\begin{proof} 
To prove that {1)$\implies$3)}, let $R\in \mathcal{L}(X_2)$ be  positive definite. Then $S:= \left[ \begin{smallmatrix}
    I & 0 \\ 0 & R
\end{smallmatrix}\right]$ is invertible and $ \tilde A_\ext=A_\ext S$. From Lemma~\ref{lem:HelpWithCompRes}, the embedding $\iota_{\mathcal{D}(A_\ext)}:\mathcal{D}(A_\ext)\to X_1\times X_2 $ is compact. Hence, the embedding $\iota_{\mathcal{D}(A_\ext S)}:\mathcal{D}(A_\ext S)\to X_1\times X_2$ factorizes as $\iota_{\mathcal{D}(A_\ext S)}= S^{-1} \iota_{\mathcal{D}(A_\ext)}S $. Thus $\iota_{\mathcal{D}(A_\ext S)}$ is compact and $\tilde A_\ext = A_\ext S$ has compact resolvent with Lemma~\ref{lem:HelpWithCompRes}.

The implication {3)$\implies$2)} is clear and {2)$\implies$1)} is similar to {1)$\implies$3)}. Thus, we have that 1), 2) and 3) are equivalent. The equivalence between 4), 5) and 7) follows from Lemma \ref{lem:CompRes} with $J=A_{12}$, whereas Lemma~\ref{lem:CompResDual} proves that 6) and 7) are equivalent. Thus, we have the equivalence of 4), 5), 6) and 7).

Now assume that 1) holds. Then $\mathcal{D}(A_\ext)$ with the graph norm is compact in $X_1\times X_2$, see Lemma~\ref{lem:HelpWithCompRes}. This implies that $\mathcal{D}(A_{12})$ and $\mathcal{D}(A_{21})$ are compact in $X_2$ and $X_1$, respectively. Then, Lemma~\ref{lem:HelpWithCompRes} implies that $\ker A_{12}$ and $\ker A_{21}$ are finite-dimensional. Thus, we have 6) and 7).

Now assume that 6) and 7) are satisfied. Then $\mathcal{D}(A_\ext) = \mathcal{D}(A_{21}) \times \mathcal{D}(A_{12})$ is compact in $X_1\times X_2$, which implies that $A_\ext$ has compact resolvent. Thus, we have 1).
\end{proof}

\section*{AI disclosure}
OpenAI ChatGPT was used to assist in constructing the counterexamples presented in Sections \ref{sec:Counterexmpl} and \ref{sec:A11=0}. The authors reviewed all the AI-assisted parts independently.




\bibliographystyle{abbrv} 
\bibliography{References}

\end{document}